\documentclass[11pt,a4paper]{article}

\usepackage{theorem,enumerate}
\usepackage{amsmath,latexsym,amssymb,amsfonts}
\usepackage{eucal}
\usepackage{color}
\usepackage{comment}
\usepackage{ulem}

\newcounter{Scounter}
\theorembodyfont{\normalfont\slshape}
\newtheorem{thm}{Theorem}[section]

\newtheorem{Thm}{Theorem}

\newtheorem{lem}[thm]{Lemma}
\newtheorem{obs}[thm]{Observation}
\newtheorem{Q}[thm]{Question} 

\newtheorem{claim}{Claim}[section]

\numberwithin{equation}{section}

\newcommand{\proof}{\medbreak\noindent\textit{Proof.}\quad}

\newcommand{\qed}{{$\quad\square$\vs{3.6}}}

\newcommand{\vs}[1]{\vspace*{#1 mm}}

\def\C{{ \mathcal{C}}}

\def\F{{ \mathcal{F}}}

\bfseries\normalfont

\title{A characterization of domination weak bicritical graphs with large diameter}

\author{
Michitaka Furuya\footnote{\texttt{michitaka.furuya@gmail.com}} \vs{5}\\
\textsl{College of Liberal Arts and Science,} \\
\textsl{Kitasato University,}\\
\textsl{1-15-1 Kitasato, Minami-ku, Sagamihara, Kanagawa 252-0373, Japan}\\
}

\date{}

\begin{document}

\maketitle

\begin{abstract}
The domination number of a graph $G$, denoted by $\gamma (G)$, is the minimum cardinality of a dominating set of $G$.
A vertex of a graph is called critical if its deletion decreases the domination number, and a graph is called critical if its all vertices are critical.
A graph $G$ is called weak bicritical if for every non-critical vertex $x\in V(G)$, $G-x$ is a critical graph with $\gamma (G-x)=\gamma (G)$.
In this paper, we characterize the connected weak bicritical graphs $G$ whose diameter is exactly $2\gamma (G)-2$.
This is a generalization of some known results concerning the diameter of graphs with a domination-criticality.
\end{abstract}



\section{Introduction}\label{sec1}

All graphs considered in this paper are finite, simple, and undirected.

Let $G$ be a graph.
We let $V(G)$ and $E(G)$ denote the vertex set and the edge set of $G$, respectively.
For $x\in V(G)$, we let $N_{G}(x)$ and $N_{G}[x]$ denote the {\it open neighborhood} and the {\it closed neighborhood} of $x$, respectively; thus $N_{G}(x)=\{y\in V(G):xy\in E(G)\}$ and $N_{G}[x]=N_{G}(x)\cup \{x\}$.
For $x,y\in V(G)$, we let $d_{G}(x,y)$ denote the {\it distance} between $x$ and $y$ in $G$.
For $x\in V(G)$ and a non-negative integer $i$, let $N^{(i)}_{G}(x)=\{y\in V(G):d_{G}(x,y)=i\}$; thus $N^{(0)}_{G}(x)=\{x\}$ and $N^{(1)}_{G}(x)=N_{G}(x)$.
The {\it diameter} of $G$, denoted by $\mbox{diam}(G)$, is defined to be the maximum of $d_{G}(x,y)$ as $x,y$ range over $V(G)$.
A vertex $x\in V(G)$ is {\it diametrical} if $\max \{d_{G}(x,y):y\in V(G)\}=\mbox{diam}(G)$.

We let $\overline{G}$ denote the {\it complement} of $G$.
For two graphs $H_{1}$ and $H_{2}$, we let $H_{1}\cup H_{2}$ denote the {\it union} of $H_{1}$ and $H_{2}$.
For a graph $H$ and a non-negative integer $s$, $sH$ denote the disjoint union of $s$ copies of $H$.
We let $K_{n}$ and $P_{n}$ denote the {\it complete graph} and the {\it path} of order $n$, respectively.

For two subsets $X,Y$ of $V(G)$, we say that $X$ {\it dominates} $Y$ if $Y\subseteq \bigcup_{x\in X}N_{G}[x]$.
A subset of $V(G)$ which dominates $V(G)$ is called a {\it dominating set} of $G$.
The minimum cardinality of a dominating set of $G$, denoted by $\gamma (G)$, is called the {\it domination number} of $G$.
A dominating set of $G$ with the cardinality $\gamma (G)$ is called a {\it $\gamma$-set} of $G$.

For terms and symbols not defined here, we refer the reader to \cite{D}.

\subsection{Motivations}\label{sec1.1}

For a given graph $G$, we can divide the set $V(G)$ into the following three subsets;
\begin{align*}
V^{0}(G)&=\{x\in V(G):\gamma (G-x)=\gamma (G)\},\\
V^{+}(G)&=\{x\in V(G):\gamma (G-x)>\gamma (G)\},\mbox{ and }\\
V^{-}(G)&=\{x\in V(G):\gamma (G-x)<\gamma (G)\}.
\end{align*}
A vertex in $V^{-}(G)$ is said to be {\it critical}.
A graph $G$ is {\it critical} if every vertex of $G$ is critical (i.e., $V(G)=V^{-}(G)$), and $G$ is {\it $k$-critical} if $G$ is critical and $\gamma (G)=k$.
Many researchers have studied critical vertices or critical graphs (for example, see~\cite{AP1,AP2,HH,S,WY}).
Among them, we focus on the following theorem which was conjectured by Brigham, Chinn and Dutton~\cite{BCD}.

\begin{Thm}[Fulman, Hanson and MacGillivray~\cite{FHM}]
\label{ThmA}
Let $k\geq 2$ be an integer, and let $G$ be a connected $k$-critical graph.
Then $\mbox{diam}(G)\leq 2k-2$.
\end{Thm}

After that, Ao~\cite{A} characterized the connected $k$-critical graphs $G$ with $\mbox{diam}(G)=2k-2$ (see Theorem~\ref{ThmE} in Subsection~\ref{sec1.2}).

Now we introduce other criticality for the domination.
A graph $G$ is {\it bicritical} if $\gamma (G-\{x,y\})<\gamma (G)$ for any pair of distinct vertices $x,y\in V(G)$, and $G$ is {\it $k$-bicritical} if $G$ is bicritical and $\gamma (G)=k$.
It is known that for $k\leq 2$, the order of a $k$-bicritical graph is at most $3$ (see~\cite{BHHR}), and hence we are interested in $k$-bicritical graphs with $k\geq 3$.
Brigham, Haynes, Henning and Rall~\cite{BHHR} gave a conjecture concerning the diameter of bicritical graphs:
For $k\geq 3$, every connected $k$-bicritical graph $G$ satisfies $\mbox{diam}(G)\leq k-1$.
However, the conjecture was disproved by the following theorem.

\begin{Thm}[Furuya~\cite{F1,F2}]
\label{ThmB}
Let $k\geq 3$ be an integer.
Then there exist infinitely many connected $k$-bicritical graphs $G$ with
$$
\mbox{diam}(G)=
\begin{cases}
3 & (k=3)\\
6 & (k=5)\\
\frac{3k-1}{2} & (k\mbox{ is odd and }k\geq 7)\\
\frac{3k-2}{2} & (k\mbox{ is even}).
\end{cases}
$$
\end{Thm}

Thus one might be interested in an upper bound of the diameter of bicritical graphs.
In~\cite{F2}, the author proved the following theorem.
(However, it is open to find a sharp upper bound of the diameter of bicritical graphs.)

\begin{Thm}[Furuya~\cite{F2}]
\label{ThmC}
Let $k\geq 3$ be an integer, and let $G$ be a connected $k$-bicritical graph.
Then $\mbox{diam}(G)\leq 2k-3$.
\end{Thm}

For convenience, let $\C$ and $\C_{B}$ denote the family of connected critical graphs and the family of connected bicritical graphs, respectively.
Here we compare Theorem~\ref{ThmA} with Theorem~\ref{ThmC}.
Although the inequalities in the theorems are similar, the two theorems are essentially different because $\C$ is different from $\C_{B}$:
\begin{enumerate}[$\bullet $]
\item
We can easily check that the graphs in $\F_{k}$ defined in Subsection~\ref{sec1.2} are critical and not bicritical.
\item
It is known that there exist infinitely many connected critical and bicritical graphs (see~\cite{BHHR,F1}), and Brigham et al.~\cite{BHHR} proved that a graph obtained from a critical and bicritical graph by expanding one vertex is bicritical and not critical.
On the other hand, there exist infinitely many connected $4$-bicritical graphs which is not critical and not obtained by the above operation (see the graph $L_{s}$ in \cite{F2}).
\end{enumerate}
In particular, $\C$ and $\C_{B}$ seems to be remotely related.

To treat the criticality and the bicriticality simultaneously, a new critical concept was defined in \cite{F2}.
A graph $G$ is {\it weak bicritical} if $V^{+}(G)=\emptyset $ and $G-x$ is critical for every $x\in V^{0}(G)$, and $G$ is {\it weak $k$-bicritical} if $G$ is weak bicritical and $\gamma (G)=k$.
Since all critical graphs and all bicritical graphs are weak bicritical, the weak bicriticality is a unification of the criticality and the bicriticality.
In~\cite{F2}, the author showed the following theorem which is a generalization of Theorem~\ref{ThmA}.

\begin{Thm}[Furuya~\cite{F2}]
\label{ThmD}
Let $k\geq 2$ be an integer, and let $G$ be a connected weak $k$-bicritical graph.
Then $\mbox{diam}(G)\leq 2k-2$.
\end{Thm}

However, Theorem~\ref{ThmC} cannot directly follow from Theorem~\ref{ThmD}.
In this paper, our main aim is to give a common generalization of Theorems~\ref{ThmA} and \ref{ThmC} by characterizing the connected weak $k$-bicritical graphs $G$ with $\mbox{diam}(G)=2k-2$.

\subsection{Main result}\label{sec1.2}

Before we state our main result, we introduce Ao's characterization.

Let $k\geq 2$ be an integer.
We define the family $\F_{k}$ of graphs as follows:
Let $m_{i}\geq 2~(1\leq i\leq k-1)$ be integers.
For each $1\leq i\leq k-1$, let $G_{i}$ be a graph isomorphic to $\overline{m_{i}K_{2}}$ (i.e., $G_{i}$ is a graph obtained from the complete graph of order $2m_{i}$ by deleting a perfect matching), and take two vertices $u_{i},v_{i}\in V(G_{i})$ with $u_{i}v_{i}\notin E(G_{i})$.
Let $G(m_{1},\ldots ,m_{k-1})$ be the graph obtained from $G_{1},\ldots ,G_{k-1}$ by identifying $v_{i}$ and $u_{i+1}$ for each $1\leq i\leq k-2$, and set
$$
\F_{k}=\{G(m_{1},\ldots ,m_{k-1}):m_{i}\geq 2,~1\leq i\leq k-1\}.
$$
By the definition of $\F_{k}$, we see the following observation.

\begin{obs}
\label{obs1.2.1}
Let $k\geq 3$, $k_{1}\geq 2$ and $k_{2}\geq 2$ be integers with $k_{1}+k_{2}-1=k$.
Then a graph $G$ belongs to $\F_{k}$ if and only if $G$ is obtained from two graphs $H_{1}\in \F_{k_{1}}$ and $H_{2}\in \F_{k_{2}}$ by identifying diametrical vertices $u_{i}$ of $H_{i}~(i\in \{1,2\})$.
\end{obs}

Ao~\cite{A} proved the following theorem.
(By using lemmas for our main result, the following theorem can be easily proved. Hence we will give its proof in Section~\ref{sec4}).

\begin{Thm}[Ao~\cite{A}]
\label{ThmE}
Let $k\geq 2$ be an integer, and let $G$ be a connected $k$-critical graph.
Then $\mbox{diam}(G)\leq 2k-2$, with the equality if and only if $G\in \F_{k}$.
\end{Thm}

Now we recursively define the family $\F^{*}_{k}~(k\geq 2)$ of graphs and the identifiable vertices of graphs in $\F^{*}_{k}$.
Let
$$
\F^{*}_{2}=\{\overline{(m+1)K_{2}},~\overline{mK_{2}\cup K_{3}},~\overline{mK_{2}\cup P_{3}}:m\geq 1\}.
$$
Note that $\F^{*}_{2}$ is equal to the family of connected weak $2$-bicritical graphs (see Lemma~\ref{lem1.22} in Subsection~\ref{sec1.3}).
For each $G\in \F^{*}_{2}$, a vertex $x\in V(G)$ is {\it identifiable} if $x\in V^{-}(G)$.
Note that if $G=\overline{(m+1)K_{2}}$, then all vertices of $G$ are identifiable; if $G=\overline{mK_{2}\cup K_{3}}$, then $G$ has exactly three non-identifiable vertices; if $G=\overline{mK_{2}\cup P_{3}}$, then $G$ has exactly two non-identifiable vertices.
We assume that $k\geq 3$, and for $2\leq k'\leq k-1$, the family $\F^{*}_{k'}$ and the identifiable vertices of graphs in $\F^{*}_{k'}$ has been defined.
Let $\F'_{k}$ be the family of graphs obtained from two graphs $H_{1}\in \F_{k_{1}}$ and $H_{2}\in \F^{*}_{k_{2}}$ with $k_{1}\geq 2$, $k_{2}\geq 2$ and $k_{1}+k_{2}-1=k$ by identifying a diametrical vertex of $H_{1}$ and an identifiable vertex of $H_{2}$.
Let $m_{i}\geq 2~(i\in \{1,2\})$, and let $u$ be the unique cut vertex of the graph $G(m_{1},m_{2})~(\in \F_{3})$.
Let $G^{1}(m_{1},m_{2})$ be the graph obtained from $G(m_{1},m_{2})$ by adding a new vertex $u'$ and joining $u'$ to all vertices in $N_{G(m_{1},m_{2})}(u)$, and let $G^{2}(m_{1},m_{2})=G^{1}(m_{1},m_{2})+uu'$.
Let
$$
\F''_{3}=\{G^{1}(m_{1},m_{2}),~G^{2}(m_{1},m_{2}):m_{i}\geq 2,~i\in \{1,2\}\},
$$
and let $\F''_{k}=\emptyset $ for $k\geq 4$.
Then by tedious argument, we see that every graph in $\F''_{3}$ is weak $3$-bicritical (but we omit detail).
Let $\F^{*}_{k}=\F'_{k}\cup \F''_{k}$ for $k\geq 3$.
For each $G\in \F^{*}_{k}$, a vertex $x\in V(G)$ is {\it identifiable} if $x\in V^{-}(G)$ and $x$ is a diametrical vertex of $G$.
By induction and Lemma~\ref{lem1.3}(ii) in Subsection~\ref{sec1.3}, we see that every graph $G\in \F^{*}_{k}$ has at least one identifiable vertex, and hence $\F^{*}_{k}$ is well-defined.
Furthermore, by the definition of $\F_{k}$ and $\F^{*}_{k}$ and Observation~\ref{obs1.2.1}, we also see that $\F_{k}\subseteq \F^{*}_{k}$ and the diameter of graphs in $\F^{*}_{k}$ is exactly $2k-2$.

Our main result is the following.

\begin{thm}
\label{thm1}
Let $k\geq 2$ be an integer, and let $G$ be a connected weak $k$-bicritical graph.
Then $\mbox{diam}(G)\leq 2k-2$, with the equality if and only if $G\in \F^{*}_{k}$.
\end{thm}

Theorem~\ref{thm1} clearly leads to Theorems~\ref{ThmA} and \ref{ThmD}.
Furthermore, it is not hard to check that no graph in $\F^{*}_{k}$ is bicritical and no graph in $\F^{*}_{k}-\F_{k}$ is critical, and so Theorem~\ref{thm1} leads to Theorems~\ref{ThmC} and \ref{ThmE}.
Therefore, Theorem~\ref{thm1} is a common generalization of some known results.

\subsection{Preliminaries}\label{sec1.3}

In this subsection, we enumerate some fundamental or preliminary results.

The following has been known property which will be used in our argument.

\begin{lem}
\label{lem1.1}
Let $G$ be a graph, and let $u,v\in V(G)$.
If $N_{G}[u]\subseteq N_{G}[v]$, then $v$ is not critical.
\end{lem}

In \cite{F2}, the author showed that the minimum degree of a connected weak bicritical graph of order at least $3$ is at least $2$.
Now we let $G$ be a disconnected weak bicritical graph.
Then we can verify that each component of $G$ is weak bicritical.
(Indeed, all components of $G$ are critical with at most one exception.)
Thus the following lemma holds.

\begin{lem}
\label{lem1.2}
Let $G$ be a weak bicritical graph, and let $G_{1}$ be a component of $G$ with $|V(G_{1})|\geq 3$.
Then the minimum degree of $G_{1}$ is at least $2$.
\end{lem}

Since the weak $1$-bicritical graphs are only $K_{1}$ and $K_{2}$, we are interested in weak $k$-bicritical graphs for $k\geq 2$.
The following lemma gives a characterization of weak $2$-bicritical graphs (or $2$-critical graphs).

\begin{lem}[Furuya~\cite{F2}]
\label{lem1.22}
A graph $G$ is weak $2$-bicritical if and only if
$$
G\in \{\overline{mK_{2}},~\overline{mK_{2}\cup K_{3}},~\overline{(m-1)K_{2}\cup P_{3}}:m\geq 1\}.
$$
In particular, a graph $G$ is $2$-critical if and only if $G\in \{\overline{mK_{2}}:m\geq 1\}$.
\end{lem}

We next focus on the coalescence of graphs.
Let $H_{1}$ and $H_{2}$ be two vertex-disjoint graphs, and let $x_{1}\in V(H_{1})$ and $x_{2}\in V(H_{2})$.
Under this notation, we let $(H_{1}\bullet H_{2})(x_{1},x_{2};x)$ denote the graph obtained from $H_{1}$ and $H_{2}$ by identifying vertices $x_{1}$ and $x_{2}$ into a vertex labeled $x$.
We call $(H_{1}\bullet H_{2})(x_{1},x_{2};x)$ the {\it coalescence} of $H_{1}$ and $H_{2}$ via $x_{1}$ and $x_{2}$.

\begin{lem}[\cite{BCD,BHHR,BS,F1}]
\label{lem1.3}
Let $H_{1}$ and $H_{2}$ be graphs, and for each $i\in \{1,2\}$, let $x_{i}$ be a non-isolated vertex of $H_{i}$.
Let $G=(H_{1}\bullet H_{2})(x_{1},x_{2};x)$.
Then the following hold.
\begin{enumerate}[{\upshape(i)}]
\item
We have $\gamma (H_{1})+\gamma (H_{2})-1\leq \gamma (G)\leq \gamma (H_{1})+\gamma (H_{2})$.
If $x_{i}$ is a critical vertex of $H_{i}$ for some $i\in \{1,2\}$, then $\gamma (G)=\gamma (H_{1})+\gamma (H_{2})-1$.
\item
If $x_{i}$ is a critical vertex of $H_{i}$ for each $i\in \{1,2\}$, then
$$
V^{-}(G)=(V^{-}(H_{1})-\{x_{1}\})\cup (V^{-}(H_{2})-\{x_{2}\})\cup \{x\}.
$$
In particular, the graph $G$ is critical if and only if both $H_{1}$ and $H_{2}$ are critical.
\end{enumerate}
\end{lem}


\section{Coalescences}\label{sec2}

In this section, we prove the following theorem.

\begin{thm}
\label{thm2.1}
Let $H_{1}$ and $H_{2}$ be graphs, and for each $i\in \{1,2\}$, let $x_{i}$ be a non-isolated vertex of $H_{i}$.
Let $G=(H_{1}\bullet H_{2})(x_{1},x_{2};x)$.
Then $G$ is weak bicritical if and only if for some $i\in \{1,2\}$,
\begin{enumerate}[{\upshape(1)}]
\item
$H_{i}$ is critical,
\item
$H_{3-i}$ is weak bicritical, and
\item
$x_{3-i}$ is a critical vertex of $H_{3-i}$.
\end{enumerate}
Furthermore, if $G$ is weak bicritical, then $\gamma (G)=\gamma (H_{1})+\gamma (H_{2})-1$.
\end{thm}
\proof
We first assume that $G$ is weak bicritical, and show that $\gamma (G)=\gamma (H_{1})+\gamma (H_{2})-1$ and (1)--(3) hold.

\begin{claim}
\label{cl2.1.1}
The vertex $x$ belongs to $V^{-}(G)$.
\end{claim}
\proof
Suppose that $x\notin V^{-}(G)$.
Then $x\in V^{0}(G)$ and $G-x$ is critical.
Since $G-x$ is the union of $H_{1}-x_{1}$ and $H_{2}-x_{2}$, $\gamma (G)=\gamma (H_{1}-x_{1})+\gamma (H_{2}-x_{2})$ and $H_{i}-x_{i}$ is critical for each $i\in \{1,2\}$.
For $i\in \{1,2\}$, let $y_{i}\in N_{H_{i}}(x_{i})$, and let $S_{i}$ be a $\gamma $-set of $H_{i}-\{x_{i},y_{i}\}$.
Then $\gamma (H_{i}-\{x_{i},y_{i}\})\leq \gamma (H_{i}-x_{i})-1$.
Since $S_{1}\cup S_{2}\cup \{x\}$ is a dominating set of $G$, we have $\gamma (H_{1}-\{x_{1},y_{1}\})+\gamma (H_{2}-\{x_{2},y_{2}\})+1=|S_{1}|+|S_{2}|+|\{x\}|\geq \gamma (G)$.
Consequently,
\begin{align*}
\gamma (G) &= \gamma (G-x)\\
&=\gamma (H_{1}-x_{1})+\gamma (H_{2}-x_{2})\\
&\geq \gamma (H_{1}-\{x_{1},y_{1}\})+\gamma (H_{2}-\{x_{2},y_{2}\})+2\\
&\geq \gamma (G)+1,
\end{align*}
which is a contradiction.
\qed

\begin{claim}
\label{cl2.1.2}
For $i\in \{1,2\}$, $x_{i}$ is a critical vertex of $H_{i}$.
\end{claim}
\proof
Let $S$ be a $\gamma $-set of $G-x$.
Then by Claim~\ref{cl2.1.1} and Lemma~\ref{lem1.3}(i), $|S|\leq \gamma (G)-1\leq \gamma (H_{1})+\gamma (H_{2})-1$.
Since $\{S\cap V(H_{1}),S\cap V(H_{2})\}$ is a partition of $S$, we have $|S\cap V(H_{i})|\leq \gamma (H_{i})-1$ for some $i\in \{1,2\}$.
Without loss of generality, we may assume that $|S\cap V(H_{1})|\leq \gamma (H_{1})-1$.
Since removing a vertex can decrease the domination number at most by one and $S\cap V(H_{1})$ is a dominating set of $H_{1}-x_{1}$, this implies that $\gamma (H_{1}-x_{1})=|S\cap V(H_{1})|=\gamma (H_{1})-1$ and $x_{1}$ is a critical vertex of $H_{1}$.
Again by Lemma~\ref{lem1.3}(i), $\gamma (G)=\gamma (H_{1})+\gamma (H_{2})-1$, and hence $|S|\leq \gamma (G)-1=\gamma (H_{1})+\gamma (H_{2})-2$.
Consequently
\begin{align*}
|S\cap V(H_{2})| &= |S|-|S\cap V(H_{1})|\\
&\leq (\gamma (H_{1})+\gamma (H_{2})-2)-(\gamma (H_{1})-1)\\
&= \gamma (H_{2})-1.
\end{align*}
Since $S\cap V(H_{2})$ is a dominating set of $H_{2}-x_{2}$, $\gamma (H_{2}-x_{2})\leq |S\cap V(H_{2})|\leq \gamma (H_{2})-1$ and $x_{2}$ is a critical vertex of $H_{2}$.
\qed

By Lemma~\ref{lem1.3} and Claim~\ref{cl2.1.2},
\begin{align}
\gamma (G)=\gamma (H_{1})+\gamma (H_{2})-1\label{eq-thm2.1-1}
\end{align}
and
\begin{align}
V^{-}(G)=(V^{-}(H_{1})-\{x_{1}\})\cup (V^{-}(H_{2})-\{x_{2}\})\cup \{x\}.\label{eq-thm2.1-2}
\end{align}
If $H_{1}$ and $H_{2}$ are critical, then (1)--(3) hold.
Thus, without loss of generality, we may assume that $H_{1}$ is not critical (i.e., $V(H_{1})-V^{-}(H_{1})\neq \emptyset $).
Let $y\in V(H_{1})-V^{-}(H_{1})$.
By (\ref{eq-thm2.1-2}), $y\notin V^{-}(G)$, and hence $G-y$ is critical.

\begin{claim}
\label{cl2.1.3}
We have $y\in V^{0}(H_{1})$.
\end{claim}
\proof
Note that $\gamma (G-\{x,y\})<\gamma (G)$, and $\gamma (H_{2}-x_{2})=\gamma (H_{2})-1$ because $x_{2}$ is a critical vertex of $H_{2}$ and removing a vertex can decrease the domination number at most by one.
Since $G-\{x,y\}$ is the union of $H_{1}-\{x_{1},y\}$ and $H_{2}-x_{2}$, this together with (\ref{eq-thm2.1-1}) leads to
\begin{align*}
\gamma (H_{1})+\gamma (H_{2})-2 &= \gamma (G)-1\\
&\geq \gamma (G-\{x,y\})\\
&= \gamma (H_{1}-\{x_{1},y\})+\gamma (H_{2}-x_{2})\\
&= \gamma (H_{1}-\{x_{1},y\})+\gamma (H_{2})-1,
\end{align*}
and so $\gamma (H_{1}-\{x_{1},y\})\leq \gamma (H_{1})-1$.
Since $S_{1}\cup \{x_{1}\}$ is a dominating set of $H_{1}-y$ for a $\gamma $-set $S_{1}$ of $H_{1}-\{x_{1},y\}$, we have
$$
\gamma (H_{1}-y)\leq \gamma (H_{1}-\{x_{1},y\})+1\leq \gamma (H_{1}).
$$
Since $y\notin V^{-}(H_{1})$, the desired conclusion holds.
\qed

Since $y$ is an arbitrary vertex in $V(H_{1})-V^{-}(H_{1})$, it suffices to show that both $H_{1}-y$ and $H_{2}$ are critical.
Note that $y\neq x_{1}$.
Now we show that
\begin{align}
\mbox{$x_{1}$ is a non-isolated vertex of $H_{1}-y$.}\label{cond-thm2.1-1}
\end{align}
By way of contradiction, we suppose that $x_{1}$ is an isolated vertex of $H_{1}-y$.
Since $x_{1}$ is a non-isolated vertex of $H_{1}$, $N_{H_{1}}(x_{1})=\{y\}$.
Since $G$ is weak bicritical and $x_{2}$ is a non-isolated vertex of $H_{2}$, the component of $G$ containing $y$ has at least three vertices.
This together with Lemma~\ref{lem1.2} implies $N_{H_{1}}(y)-\{x_{1}\}\neq \emptyset $.
Let $y'\in N_{H_{1}}(y)-\{x_{1}\}$.
Since $G-y$ is critical, $\gamma (G-\{y,y'\})\leq \gamma (G)-1=\gamma (H_{1})+\gamma (H_{2})-2$.
Let $S$ be a $\gamma $-set of $G-\{y,y'\}$.
If $x\in S$, let $S'=((S-\{x\})\cap V(H_{2}))\cup \{x_{2}\}$; if $x\notin S$, let $S'=S\cap V(H_{2})$.
In either case, $S'$ is a dominating set of $H_{2}$, and hence $|(S-\{x\})\cap V(H_{1})|=|S|-|S'|\leq (\gamma (H_{1})+\gamma (H_{2})-2)-\gamma (H_{2})=\gamma (H_{1})-2$.
Since $(S-\{x\})\cap V(H_{1})$ is a dominating set of $H_{1}-\{x,y,y'\}$, $S''=((S-\{x\})\cap V(H_{1}))\cup \{y\}$ is a dominating set of $H_{1}$ with $|S''|\leq \gamma (H_{1})-1$, which is a contradiction.
Thus (\ref{cond-thm2.1-1}) holds.

Recall that $G-y$ is critical.
Since $G-y=((H_{1}-y)\bullet H_{2})(x_{1},x_{2};x)$, it follows from Lemma~\ref{lem1.3}(ii) and (\ref{cond-thm2.1-1}) that $H_{1}-y$ and $H_{2}$ are critical.

We next assume that (1)--(3) hold, and show that $G$ is weak bicritical.
We may assume that $i=1$ (i.e., $H_{1}$ is critical, $H_{2}$ is weak bicritical, and $x_{2}$ is a critical vertex of $H_{2}$).
By Lemma~\ref{lem1.3}(i), $\gamma (G)=\gamma (H_{1})+\gamma (H_{2})-1$.
If $G$ is critical, then the desired conclusion holds.
Thus $V(G)-V^{-}(G)\neq \emptyset $.
Let $y\in V(G)-V^{-}(G)$.
By Lemma~\ref{lem1.3}(ii), $y\in V^{0}(H_{2})$, and hence $H_{2}-y$ is critical.

\begin{claim}
\label{cl2.1.4}
We have $y\in V^{0}(G)$.
\end{claim}
\proof
Let $S_{1}$ be a $\gamma $-set of $H_{1}$, and let $S_{2}$ be a $\gamma $-set of $H_{2}-\{x_{2},y\}$.
If $x_{1}\in S_{1}$, let $S=(S_{1}-\{x_{1}\})\cup S_{2}\cup \{x\}$; if $x_{1}\not\in S_{1}$, let $S=S_{1}\cup S_{2}$.
In either case, $S$ is a dominating set of $G-y$.
Since $|S|=\gamma (H_{1})+\gamma (H_{2}-\{x_{2},y\})\leq \gamma (H_{1})+(\gamma (H_{2})-1)=\gamma (G)$, we have $\gamma (G-y)\leq \gamma (G)$.
Since $y\notin V^{-}(G)$, the desired conclusion holds.
\qed

Since $y$ is an arbitrary vertex in $V(G)-V^{-}(G)$, it suffices to show that $G-y$ is critical.
Note that $y\neq x$.
Now we show that
\begin{align}
\mbox{$x_{2}$ is a non-isolated vertex of $H_{2}-y$.}\label{cond-thm2.1-2}
\end{align}
Recall that $x_{2}$ is a non-isolated vertex of $H_{2}$.
Furthermore, since $x_{2}$ is a critical vertex of $H_{2}$, the component of $H_{2}$ containing $x_{2}$ is not isomorphic to $K_{2}$, and hence the component of $H_{2}$ containing $x_{2}$ has at least three vertices.
This together with Lemma~\ref{lem1.2} implies that the degree of $x_{2}$ in $H_{2}$ is at least $2$, and so the degree of $x_{2}$ in $H_{2}-y$ is at least $2$.
Thus (\ref{cond-thm2.1-2}) holds.

Recall that both $H_{1}$ and $H_{2}-y$ are critical.
Since $G-y=(H_{1}\bullet (H_{2}-y))(x_{1},x_{2};x)$, it follows from Lemma~\ref{lem1.3}(ii) and (\ref{cond-thm2.1-2}) that $G-y$ is critical.

This completes the proof of Theorem~\ref{thm2.1}.
\qed

\section{Sufficient pairs}\label{sec3}

Let $l\geq 3$ be an integer, and let $G$ be a connected graph.
A pair $(x,j)$ of a vertex $x\in V(G)$ and an integer $j\geq 2$ is {\it $l$-sufficient} if $x$ is a diametrical vertex of $G$ and there exists a $\gamma $-set $S$ of $G$ with $|S\cap (\bigcup _{0\leq i\leq j}N^{(i)}_{G}(x))|\geq (j+l)/2$.

\begin{lem}[Furuya~\cite{F2}]
\label{Lem3-A}
Let $k\geq 3$ and $l\geq 3$ be integers, and let $G$ be a connected weak $k$-bicritical graph having an $l$-sufficient pair.
Then $\mbox{diam}(G)\leq 2k-l+1$.
\end{lem}

\begin{thm}
\label{thm3.1}
Let $k\geq 3$ be an integer, and let $G$ be a connected weak $k$-bicritical graph.
If $G$ has a diametrical vertex $x$ such that $\bigcup _{1\leq i\leq 3}N^{(i)}_{G}(x)\subseteq V^{-}(G)$ and $|N^{(2)}_{G}(x)|\geq 2$, then $\mbox{diam}(G)\leq 2k-3$.
\end{thm}
\proof
We show that $\mbox{diam}(G)\leq 3$ or $G$ has a $4$-sufficient pair.
By way of contradiction, we suppose that $\mbox{diam}(G)\geq 4$ and $G$ has no $4$-sufficient pair.
For each $i\geq 0$, let $X_{i}=N^{(i)}_{G}(x)$ and $U_{i}=X_{0}\cup X_{1}\cup \cdots \cup X_{i}$.

\begin{claim}
\label{cl3.1.1}
If a set $S\subseteq V(G)$ dominates $N_{G}[x]$ and $|S\cap U_{2}|\leq 1$, then $x$ is the unique vertex of $S\cap U_{2}$.
\end{claim}
\proof
By the assumption of the claim, there exists a vertex $z\in N_{G}[x]$ dominating $N_{G}[x]$ in $G$.
Since $N_{G}[x]\subseteq N_{G}[z]$, if $z\neq x$, then $z\in N^{(1)}_{G}(x)$ and $z$ is not a critical vertex of $G$ by Lemma~\ref{lem1.1}, which contradicts the assumption of the theorem.
\qed

Let $w_{2},w'_{2}\in X_{2}$ be distinct vertices, and let $S_{1}$ be a $\gamma $-set of $G-w_{2}$.
Note that $S_{1}\cup \{w_{2}\}$ is a $\gamma $-set of $G$.
Since $G$ has no $4$-sufficient pair, $|(S_{1}\cup \{w_{2}\})\cap U_{2}|<(2+4)/2=3$, and so $|S_{1}\cap U_{2}|\leq 1$.
Since $S_{1}$ dominates $N_{G}[x]$ in $G$, it follows from Claim~\ref{cl3.1.1} that $x$ is the unique vertex in $S_{1}\cap U_{2}$.
Since $G$ has no $4$-sufficient pair, $|(S_{1}\cup \{w_{2}\})\cap U_{4}|<(4+4)/2=4$, and so $|S_{1}\cap U_{4}|\leq 2$.
Since $|X_{2}|\geq 2$ and $S_{2}$ dominates $(X_{2}\cup X_{3})-\{w_{2}\}$, there exists a vertex $w_{3}\in X_{3}$ dominating $(X_{2}\cup X_{3})-\{w_{2}\}$ in $G-w_{2}$.

Let $S_{2}$ be a $\gamma $-set of $G-w_{3}$.
Note that $S_{2}\cup \{w'_{2}\}$ is a $\gamma $-set of $G$ because $w_{3}w'_{2}\in E(G)$.
Since $G$ has no $4$-sufficient pair, $|(S_{2}\cup \{w'_{2}\})\cap U_{2}|<(2+4)/2=3$, and so $|S_{2}\cap U_{2}|\leq 1$.
Since $S_{2}$ dominates $N_{G}[x]$ in $G$, it follows from Claim~\ref{cl3.1.1} that $x$ is the unique vertex in $S_{2}\cap U_{2}$.
Since $G$ has no $4$-sufficient pair, $|(S_{2}\cup \{w'_{2}\})\cap U_{4}|<(4+4)/2=4$, and so $|S_{2}\cap U_{4}|\leq 2$.
Since $S_{2}$ dominates $(X_{2}\cup X_{3})-\{w_{3}\}$, there exists a vertex $w'_{3}\in X_{3}$ dominating $(X_{2}\cup X_{3})-\{w_{3}\}$ in $G-w_{3}$.
Recall that $w_{3}$ dominates $X_{3}$ in $G-w_{2}$.
Thus $w_{3}w'_{3}\in E(G)$, and hence $S_{2}$ is a dominating set of $G$, which is a contradiction.

Consequently $\mbox{diam}(G)\leq 3$ or $G$ has a $4$-sufficient pair.
In either case, it follows from Lemma~\ref{Lem3-A} that the desired conclusion holds.
\qed

\section{Proof of Theorems~\ref{ThmE} and \ref{thm1}}\label{sec4}

In this section, we prove Theorems~\ref{ThmE} and \ref{thm1}.
As we mentioned in Subsection~\ref{sec1.2}, $\F_{k}\subseteq \F^{*}_{k}$ and the diameter of graphs in $\F^{*}_{k}$ is exactly $2k-2$.
By Lemma~\ref{lem1.22}, $\F_{2}$ is equal to the family of connected $2$-critical graphs.
Thus by induction and Lemma~\ref{lem1.3}(ii), we see that all graphs in $\F_{k}$ are $k$-critical, and so
\begin{align}
\mbox{if a graph $G$ belongs to $\F_{k}$, then $G$ is $k$-critical and $\mbox{diam}(G)=2k-2$.}\label{eq-sec4-1}
\end{align}
Recall that every graph in $\F^{*}_{2}$ is weak $2$-bicritical and every graph in $\F''_{3}$ is weak $3$-bicritical.
This together with induction and Theorem~\ref{thm2.1} implies that all graphs in $\F^{*}_{k}$ are weak $k$-bicritical, and so
\begin{align}
\mbox{if a graph $G$ belongs to $\F^{*}_{k}$, then $G$ is weak $k$-bicritical and $\mbox{diam}(G)=2k-2$.}\label{eq-sec4-2}
\end{align}

\medbreak\noindent\textit{Proof of Theorem~\ref{ThmE}.}\quad
Let $k$ and $G$ be as in Theorem~\ref{ThmE}.
By (\ref{eq-sec4-1}), it suffices to show that
\begin{align}
\mbox{if $\mbox{diam}(G)\geq 2k-2$, then $G\in \F_{k}$.}\label{ThmE-1}
\end{align}
We proceed by induction on $k$.

If $k=2$, then Lemma~\ref{lem1.22} leads to (\ref{ThmE-1}).
Thus we may assume that $k\geq 3$.
Suppose that $\mbox{diam}(G)\geq 2k-2$.
Let $w$ be a diametrical vertex of $G$.
If $|N^{(2)}_{G}(w)|\geq 2$, then $\mbox{diam}(G)\leq 2k-3$ by Theorem~\ref{thm3.1}, which is a contradiction.
Thus $|N^{(2)}_{G}(w)|=1$.
In particular, $G$ has a cut vertex $x$.
Hence we can write $G$ as $G=(H_{1}\bullet H_{2})(x_{1},x_{2};x)$ for two graphs $H_{1}$ and $H_{2}$ and vertices $x_{i}\in V(H_{i})~(i\in \{1,2\})$.
For each $i\in \{1,2\}$, set $k_{i}=\gamma (H_{i})$.
By Lemma~\ref{lem1.3}, $H_{1}$ and $H_{2}$ are critical and $k_{1}+k_{2}-1=\gamma (H_{1})+\gamma (H_{2})-1=\gamma (G)=k$.
Furthermore, we have $\mbox{diam}(G)\leq \mbox{diam}(H_{1})+\mbox{diam}(H_{2})$.
By induction hypothesis, $\mbox{diam}(H_{i})\leq 2k_{i}-2$, with the equality if and only if $H_{i}\in \F_{k_{i}}$.
Consequently, we have $2k-2\leq \mbox{diam}(G)\leq (2k_{1}-2)+(2k_{2}-2)=2k-2$.
This implies that $H_{i}\in \F_{k_{i}}$ and $x_{i}$ is a diametrical vertex of $H_{i}$.
Then by Observation~\ref{obs1.2.1}, we have $G\in \F_{k}$.

This completes the proof of Theorem~\ref{ThmE}.
\qed

\medbreak\noindent\textit{Proof of Theorem~\ref{thm1}.}\quad
Let $k$ and $G$ be as in Theorem~\ref{thm1}.
By (\ref{eq-sec4-2}), it suffices to show that
\begin{align}
\mbox{if $\mbox{diam}(G)\geq 2k-2$, then $G\in \F^{*}_{k}$.}\label{thm1-1}
\end{align}
We proceed by induction on $k$.

If $k=2$, then Lemma~\ref{lem1.22} leads to (\ref{thm1-1}).
Thus we may assume that $k\geq 3$.
Suppose that $\mbox{diam}(G)\geq 2k-2$.
If $G$ is critical, then it follows from Theorem~\ref{ThmE} that $G\in \F_{k}~(\subseteq \F^{*}_{k})$, as desired.
Thus we may assume that $G$ is not critical (i.e., $V^{0}(G)\neq \emptyset $).
Let $w,w'\in V(G)$ be vertices with $d_{G}(w,w')=\mbox{diam}(G)$.

\begin{claim}
\label{cl4.1}
If $G$ has no cut vertex, then $G\in \F^{*}_{k}$.
\end{claim}
\proof
Note that $|N^{(2)}(w)|\geq 2$.
If $V^{0}(G)\subseteq \{w,w'\}$ (i.e., $V(G)-\{w,w'\}\subseteq V^{-}(G)$), then by Theorem~\ref{thm3.1}, we have $\mbox{diam}(G)\leq 2k-3$, which is a contradiction.
Thus $V^{0}(G)-\{w,w'\}\neq \emptyset $.
Let $z\in V^{0}(G)-\{w,w'\}$.
Then $G-z$ is a connected critical graph and
$$
\mbox{diam}(G-z)\geq d_{G-z}(w,w')\geq d_{G}(w,w')=\mbox{diam}(G)\geq 2k-2.
$$
This together with Theorem~\ref{ThmE} forces $G-z\in \F_{k}$ and $\mbox{diam}(G-z)=d_{G-z}(w,w')=\mbox{diam}(G)=2k-2$.
By the definition of $\F_{k}$, we have $|N^{(2)}_{G-z}(w)|=|N^{(4)}_{G-z}(w)|=1$.
Write $N^{(2)}_{G-z}(w)=\{z'\}$.
Since $G$ has no cut vertex, the following hold:
\begin{enumerate}[$\bullet $]
\item
$k=3$,
\item
$z$ is adjacent to a vertex in $N^{(1)}_{G-z}(w)$ and a vertex in $N^{(3)}_{G}(w)$, and
\item
$N_{G}(z)\subseteq \bigcup _{1\leq i\leq 3}N^{(i)}_{G-z}(w)$.
\end{enumerate}

Suppose that $z'$ is a critical vertex of $G$, and let $S$ be a $\gamma $-set of $G-z'$.
Since $N_{G}(z)\subseteq N_{G}[z']$ and $S$ is not a dominating set of $G$, this forces $zz'\notin E(G)$ and $z\in S$.
Since $S$ dominates $w$, $S\cap N_{G}[w]\neq \emptyset $.
In particular, $|(S\cup \{z'\})\cap (\bigcup _{0\leq i\leq 2}N^{(i)}_{G}(w))|\geq 3$.
Since $S\cup \{z'\}$ is a $\gamma $-set, $(w,2)$ is a $4$-sufficient pair.
This together with Lemma~\ref{Lem3-A} implies that $\mbox{diam}(G)\leq 2k-3$, which is a contradiction.
Thus $z'$ is not a critical vertex of $G$ (i.e., $z'\in V^{0}(G)$).

Replacing the role of $z$ and $z'$, we have $G-z'\in \F_{k}$ and $N_{G-z'}(z)=N^{(1)}_{G-z'}(w)\cup N^{(3)}_{G-z'}(w)$.
Hence $G$ is isomorphic to a graph in $\F''_{3}~(\subseteq \F^{*}_{3})$.
\qed

By Claim~\ref{cl4.1}, we may assume that $G$ has a cut vertex $x$.
Then we can write $G$ as $G=(H_{1}\bullet H_{2})(x_{1},x_{2};x)$ for two graphs $H_{1}$ and $H_{2}$ and vertices $x_{i}\in V(H_{i})~(i\in \{1,2\})$.
For each $i\in \{1,2\}$, set $k_{i}=\gamma (H_{i})$.
Having Theorem~\ref{thm2.1} in mind, we may assume that $H_{1}$ is critical, $H_{2}$ is weak bicritical and $x_{2}$ is a critical vertex of $H_{2}$.
Furthermore, $k_{1}+k_{2}-1=\gamma (H_{1})+\gamma (H_{2})-1=\gamma (G)=k$.
By induction hypothesis, $\mbox{diam}(H_{1})\leq 2k_{1}-2$, with the equality if and only if $H_{1}\in \F_{k_{1}}$.
By Theorem~\ref{ThmE}, $\mbox{diam}(H_{2})\leq 2k_{2}-2$, with the equality if and only if $H_{2}\in \F^{*}_{k_{2}}$.
Since $\mbox{diam}(G)\leq \mbox{diam}(H_{1})+\mbox{diam}(H_{2})$, we have $2k-2\leq \mbox{diam}(G)\leq (2k_{1}-2)+(2k_{2}-2)=2k-2$.
This implies that $H_{1}\in \F_{k_{1}}$, $H_{2}\in \F^{*}_{k_{2}}$ and $x_{i}$ is a diametrical vertex of $H_{i}$.
Since $x_{2}$ is a critical vertex of $H_{2}$, it follows from the definition of $\F^{*}_{k}$, we have $G\in \F^{*}_{k}$.

This completes the proof of Theorem~\ref{thm1}.
\qed

\section*{Acknowledgment}
This work was supported by JSPS KAKENHI Grant number 26800086.


\begin{thebibliography}{99}
\bibitem{AP1}
N.~Ananchuen and M.D.~Plummer,
Matchings in $3$-vertex-critical graphs: the even case,
{\it Networks} {\bf 45} (2005) 210--213.

\bibitem{AP2}
N.~Ananchuen and M.D.~Plummer,
Matchings in $3$-vertex-critical graphs: the odd case,
{\it Discrete Math.} {\bf 307} (2007) 1651--1658.

\bibitem{A}
S.~Ao,
Independent domination critical graphs,
Masters Thesis, University of Victoria, Victoria, BC, Canada, 1994.

\bibitem{BCD}
R.C.~Brigham, P.Z.~Chinn and R.D.~Dutton,
Vertex domination-critical graphs,
{\it Networks} {\bf 18} (1988) 173--179.

\bibitem{BHHR}
R.C.~Brigham, T.W.~Haynes, M.A.~Henning and D.F.~Rall,
Bicritical domination,
{\it Discrete Math.} {\bf 305} (2005) 18--32.

\bibitem{BS}
T.~Burton and D.P.~Sumner,
Domination dot-critical graphs,
{\it Discrete Math.} {\bf 306} (2006) 11--18.

\bibitem {D}
R.~Diestel,
{\it Graph Theory (4th edition)},
Graduate Texts in Mathematics \textbf{173}, Springer, 2010.

\bibitem{FHM}
J.~Fulman, D.~Hanson and G.~MacGillivray,
Vertex domination-critical graphs,
{\it Networks} {\bf 25} (1995) 41--43.

\bibitem{F1}
M.~Furuya,
Construction of $(\gamma ,k)$-critical graphs,
Australas. J. Combin. {\bf 53} (2012) 53--65.

\bibitem{F2}
M.~Furuya,
On the diameter of domination bicritical graphs,
{\it Australas. J. Combin.} {\bf 62} (2015) 184--196.

\bibitem{HH}
T.W.~Haynes and M.A.~Henning,
Changing and unchanging domination: a classification,
{\it Discrete Math.} {\bf 272} (2003) 65--79.

\bibitem{S}
V.~Samodivkin,
Changing and unchanging of the domination number of a graph,
{\it Discrete Math.} {\bf 308} (2008) 5015--5025.

\bibitem{WY}
T.~Wang and Q.~Yu,
A conjecture on $k$-factor-critical and $3$-critical graphs,
{\it Sci. China Math.} {\bf 53} (2010) 1385--1391.

\end{thebibliography}
\end{document}